\documentclass{amsart}

\newtheorem{theorem}{Theorem}[section]

\newtheorem{proposition}[theorem]{Proposition}

\theoremstyle{definition}
\newtheorem{definition}[theorem]{Definition}
\newtheorem{example}[theorem]{Example}

\theoremstyle{remark}

\numberwithin{equation}{section}

\begin{document}

\title{ON MULTIPLE LEBESGUE FUNCTIONS}

\author{Oleg\,N.~ Ageev
 }
\address{Department of Mathematics and Mechanics, Lomonosov Moscow State
University, Leninskiye Gory, Main Building, GSP-1, 119991 Moscow,
Russia } \email{AgeevOlegs@gmail.com}
% \curraddr{Max Planck Institute of Mathematics,
% P.O.Box 7280, D-53072 Bonn, Germany} \email{ageev@mpim-bonn.mpg.de}
\thanks{}

%    General info

\date{February 14. 2017.}

%\dedicatory{(Communicated by ?????)}

\keywords{Multiple mixing, Lebesgue spectrum.}

\begin{abstract}
We introduce a notion being a $k$-fold Lebesgue function for measure
preserving transformations, where any $2$-fold Lebesgue function is
just ordinary Lebesgue. We discuss how this new metrical
isomorphisms invariant of dynamical systems is related to others
classical notions in ergodic theory, mostly focusing on its spectral
aspects. In particular, for transformations with sufficiently many
multiple Lebesgue functions we treat the corresponding multiple
analogs of very well-known problems of Banach and Rokhlin.
\end{abstract}

\maketitle

\section*{}
\section{Introduction}

By a transformation $T$ we mean an invertible measure-preserving map
defined on a non-atomic standard Borel probability space
$(X,\mathcal{F} , \mu)$.  The \textbf{spectral} properties of $T$
are those of the induced ( Koopman's) unitary operator on
 $L^2( \mu)$ defined by
\[
\widehat{T}:L^2( \mu)\rightarrow L^2( \mu);\quad
\widehat{T}f(x)=f(Tx).
\]
In the sequel, we keep the notation $L_0^2( \mu)$ for the main
$\widehat{T}$-invariant subspace of \textbf{centered } functions $f$
(i.e. $\mu (f)=0$) in $L^2( \mu)$ and $C(f_1,f_2,\ldots )$, called
\textbf{ the component }, for the minimal closed
$\widehat{T}$-invariant subspace in $L^2( \mu)$, containing
functions $f_1,f_2,\ldots$. Usually, in ergodic theory, we talk
about reduced spectral invariants of $\widehat{T}$ that are
restricted to $L_0^2( \mu)$ because of the trivial rest. For
example, the oldest unsolved problem, associated with Banach, asks
whether there exists a transformation with simple Lebesgue spectrum
(i.e. whether there exists a Lebesgue function $f$ such that
$C(f)=L_0^2( \mu)$). Let us remind that by the Bochner theorem we
can find to each $\phi \in L^2( \mu)$ a unique (spectral) measure
$\sigma_\phi$ on the unique circle of the complex plain
$\mathbb{T}=\{\lambda \in \mathbb{C} : |\lambda|=1\}$ with the
following moments.
\[
\sigma_\phi (n)=\int_\mathbb{T}\lambda^nd\sigma_\phi(\lambda)=\int_X
\widehat{T}^n\phi\cdot\overline{\phi}d\mu,  \ n \in \mathbb{Z}.
\]
And, as above, in the classical setting, $f$ to be a Lebesgue
function
 means just $\sigma_f$ to be the Lebesgue measure or,
equivalently,
\[
\int_X \widehat{T}^nf\cdot\overline{f}d\mu=0
\]
if $n\neq 0$, and $\int_X |f|^2d\mu=2\pi$.

It is well known that any Koopman operator is just a unitary
operator that is additionally  multiplicative on the set of bounded
functions. Besides, needless to say that it is still unknown which
unitary operator can be realized as Koopman's one. There were
invented many other invariants of transformations, but the
classification problem is unsolved even in certain subclasses as,
for example, the sets of rank one or mixing transformations. It all
can be considered as an easy motivation in study of other
invariants.

One of the ways to get involved some $piece$ of the multiplicative
structure is to care about higher order moments. It is not a new
isomorphisms invariant, but still not in active use, because,
probably, of harder calculation.

One more our motivation is to look a bit differently at the class of
transformations with Lebesgue spectrum keeping in mind the nice Host
result (see [9]) that if there is a counterexample to Rokhlin's
problem on multiple mixing, then it must have at least an absolutely
continuous component in its spectrum.

To understand the motivation better, let us rewrite classical
notions mixing and multiple mixing according to the functional point
of view.

\subsection{Mixing case}
From now on we keep for mixing the following version of the
definition. We say that a transformation $T$ is mixing iff for any
  $f,g \in L^2( \mu)$
\[
\int_X \widehat{T}^nf\cdot\overline{g}d\mu \rightarrow \int_X
fd\mu\cdot\int_X \overline{g}d\mu \ \mbox{as } n\rightarrow \infty,
\eqno (1)
\]
or, in an equivalent form, for any
  $f,g \in L_0^2( \mu)$
\[
\int_X \widehat{T}^nf\cdot\overline{g}d\mu \rightarrow 0 \ \mbox{as
} n\rightarrow \infty.
\]

Following, for example, [3], and ..., we say that $T$ is mixing in
some component $B=C(f_1,f_2,\ldots)$ (and $B$ is said to be a
\textbf{mixing component} ) iff for any
  $f,g \in B$ the assumption (1) is true.

  For every transformation $T$ there exist many
  decompositions of the space $L_0^2( \mu)$ into the orthogonal sum
  of at most countably many components of the form $C(f_i)$
  ($i=1,2,\ldots$) that are called
\textbf{cyclic}. Obviously, $T$ is mixing if  all the components
$C(f_i)$ are mixing ones. This implies that to check being mixing
for $T$ is just to check that the moments $\sigma_{f_i} (n)$ vanish
for any $i$.

The notion being $k$-fold mixing ($k\geq2$) for a transformation $T$
can be defined as follows. For any collection of $k$ bounded
functions $f_1,f_2,\ldots,f_k$ in $L^2( \mu)$ and for any sequence
of $k-1$ integers $n_1(m),n_2(m),\ldots,n_{k-1}(m)$
($m=1,2,\ldots$), if $n_i(m)-n_j(m)\rightarrow\infty$ and
$n_i(m)\rightarrow\infty$ as $m\rightarrow +\infty$  for any $i\neq
j$, then
\[
\int_X\widehat{T}^{n_1(m)}f_1\cdot \widehat{T}^{n_2(m)}f_2 \cdots
\widehat{T}^{n_{k-1}(m)}f_{k-1} \cdot \overline{f}_k d\mu
\rightarrow \eqno (2)
\]
\[
\int_X f_1d\mu\cdot\int_X f_2d\mu\cdots \int_X
f_{k-1}d\mu\cdot\int_X \overline{f}_kd\mu \ \mbox{as } m\rightarrow
+\infty.
\]

Obviously, $k$-fold mixing implies mixing of all smaller orders, and
the famous unsolved Rokhlin problem on multiple mixing asks whether
there exists a ($2$-fold) mixing transformation which is not mixing
of all orders (i.e. not multiple mixing).

\begin{definition} Following [14], we say that a transformation $T$
is $k$-fold mixing ($k\geq 2$)  \textbf{ functionally} in some
component $C(g_1,g_2,\ldots)$ if there exists a dense subset $A$ of
$C(g_1,g_2,\ldots)$ such that (2) is satisfied for any collection of
$k$ functions $f_1,f_2,\ldots,f_k$ in $A$.

\end{definition}
Trivially, all the above definitions are equivalent for $k=2$.
Nevertheless, for higher orders it is a tricky question, because
there is not much known here. To avoid such unpleasant things as,
say, treating of questions how many bounded functions some component
may have and what is the set of functions in a fixed component such
that integrals in (2) exist for sufficiently large $m$, we slightly
modify the set $A$\footnote{ In fact,in [14], there was no taking
into account that integrals in (2) may not exist for some functions.
Therefore, we really need some $A$.}.

\begin{definition}  We say that a transformation $T$
is $k$-fold mixing ($k\geq 2$)  in some component
$C(g_1,g_2,\ldots)$ if it is $k$-fold mixing functionally for the
subset $A$, where $A$ contains any bounded function  in
$C(g_1,g_2,\ldots)$.

\end{definition}
One can easily construct transformations (see Example 2.2 below)
with some orthogonal decomposition of $ L^2_0( \mu)$ into $3$-fold
mixing cyclic components $C(f_i)$, where every $C(f_i)$ is not
mixing. Thus the transformations are not $3$-fold mixing. Besides,
in spite of the case $k=2$,  $3$-fold mixing even in $ L^2_0( \mu)$
does not imply $3$-fold mixing in  the whole $ L^2( \mu)$ (see
Example 2.3 ). Let us also mention that there exist transformations
with  a mixing component that is not $3$-fold mixing and contains a
dense set of bounded functions ( see [14]). Adding all we see that
there is not much of the expected similarities  for higher orders of
mixing if we do not add some regularities conditions between
components and the constant. One of the ways to do it is the
following.

Fix a transformation $T$, and $k\geq 2$. Take some orthogonal
decomposition
\[
L^2_0( \mu)=\bigoplus_{i\geq 1}C(f_i), \ \mbox{then } L^2(
\mu)=\bigoplus_{i\geq 0}C(f_i),
\]
where $f_0$ is a non-zero constant function. Let us say that this
decomposition is  \textbf{sufficiently ($k$-)good} if

\[
\forall i \int_X |f_i|^{2k-2}d\mu<\infty.
\]

For any sufficiently good decomposition , it can easily be checked
that the transformation $T$ is $k$-fold mixing if and only if for
any collection of indexes $i_1,i_2,\ldots,i_k$
($i^2_1+i^2_2+\ldots+i^2_k\neq 0$) and for any sequence of $k-1$
integers $n_1(m),n_2(m),\ldots,n_{k-1}(m)$ ($m=1,2,\ldots$), if
$n_i(m)-n_j(m)\rightarrow\infty$ and $n_i(m)\rightarrow\infty$ as
$m\rightarrow +\infty$  for any $i\neq j$, then
\[
\int_X\widehat{T}^{n_1(m)}f_{i_1}\cdot \widehat{T}^{n_2(m)}f_{i_2}
\cdots \widehat{T}^{n_{k-1}(m)}f_{i_{k-1}} \cdot \overline{f}_{i_k}
d\mu \rightarrow 0 \ \mbox{as } m\rightarrow +\infty.\eqno (3)
\]

Similarly, take a transformation $T$ and its cyclic component
$C(f)\subseteq L^2_0( \mu)$. Suppose $\int_X |f|^{2k-2}d\mu<\infty$
and $C(f)$ contains a dense subset of bounded functions. Then, $T$
is $k$-fold mixing in the component $C(f)$ if and only if (3) are
all true, where we replaced every $f_{i_j}$ ($j=1,2,\ldots,k$) by
$f$.

 By the spectral
theorem, for any transformation with simple spectrum, we can write
$L^2_0( \mu)=C(f)$ for some $f$. Applying the Alexeyev theorem (see
[5]), $f$ can be taken bounded. In the general case, we stress
 that there is still no examples of (even not necessarily $k$-fold
mixing) transformations without sufficiently good decompositions of
$L^2_0( \mu)$.

\subsection{Lebesgue case}

 \begin{definition} We say that a non-zero function $f \in L^2( \mu)$
 is
 \textbf{ $k$-fold Lebesgue} ($k\geq 2$) if for any collection
  ($n_1,n_2,\ldots,n_{k-1}$) of pairwise different non-zero
integers and for any $m\leq k$
 \begin{description}

  \item [(i)] \  \ $\widehat{T}^{n_1}f\cdot \widehat{T}^{n_2}f \cdots \widehat{T}^{n_{m-1}}f \in L^2(
  \mu)$,
  \\
  \item [(ii)] \  \  $\int_X\widehat{T}^{n_1}f\cdot \widehat{T}^{n_2}f
   \cdots \widehat{T}^{n_{k-1}}f \cdot \overline{f} d\mu=0$.
\end{description}

 \end{definition}

 Obviously, for $k=2$ the only difference is we admit to get any positive
 $L^2$-norm for Lebesgue functions.
 We remark that in the theory of real-valued random processes
 the above assumption (i)  is somehow related to the most common
 finiteness restriction of the $k$$\frac{th}{}$ moment (i.e.
 $\int_X |f|^kd\mu < \infty$). Thus developing of the
 corresponding theory is a bit of different and independent
 interest. However, in our opinion, (i) is most adapted to the
 purposes of this paper\footnote{  By spectral arguments  we can prove
 that for any transformation holding at least one Lebesgue function there is
 a lot of real-valued different ones, but we do not believe that the analogous statement
 is true even for $3$-fold Lebesgue functions. }.

 We begin with the following simple observation.
\begin{proposition} Let $T$ be any transformation and $f$ be any
 $k$-fold Lebesgue function
for $T$. The following are true:
\begin{description}

  \item [(i)] \  \ The function $\overline{f}$ is  $k$-fold Lebesgue as
  well.
  \\
  \item [(ii)] \  \ If \ $T$ is ergodic, then $f \in L_0^2( \mu)$.
  \\
  \item [(iii)] \  \ If $f \in L_0^2( \mu)$, then
  $T$  is  $k$-fold mixing functionally in $C(f)$ .
  \end{description}
\end{proposition}

In the sequal, let us restrict ourselves  to the first non-trivial
case $k=3$ mostly. We do hope all we deal below with can be
naturally extended to the general case.

Now we introduce a notion of some $k$-fold analog of transformations
with Lebesgue spectrum.

According to the mixing case discussion above, it looks natural to
leave the following $naive$ definition.

Fix some collection $\{f_i\}_{i\in I}$, ($I\subseteq \mathbb{N} $)
such that $L_0^2( \mu)=C(\{f_i\}_{i\in I})$ for a transformation
$T$, and every $f_i$ is a $3$-fold Lebesgue function. Denote
$B=\{\widehat{T}^nf_i : i\in I, n\in \mathbb{Z}\}\cup \{f_0\} $,
where $f_0$ is a fixed non-zero constant function, and any element
$\widehat{T}^{n_1}f_{i_1}$ is different to
$\widehat{T}^{n_2}f_{i_2}$ iff $n_1\neq n_2$ or $i_1\neq i_2$ (they
can be equal as elements of $L_0^2( \mu)$). One can say that the
transformation $T$ has \textbf{$3$-fold Lebesgue functions generated
spectrum } if for any collection $\varphi_i,\varphi_j,\varphi_k$ of
pairwise different elements of $B$
\[
\int_X\varphi_i\cdot\varphi_j\cdot\overline{\varphi_k}d\mu=0.
\]

At the first side this definition looks attractive, since, for
example, $3$-fold Lebesgue functions generated spectrum implies
$2$-fold Lebesgue functions generated spectrum. The first reason why
this definition is not so good is the fact that for $k=2$ it gives
the homogeneous Lebesgue spectrum only (i.e. $L_0^2( \mu)$ is the
orthogonal sum of cyclic Lebesgue components). Besides, unitary
operators with Lebesgue spectrum (i.e. with Lebesgue measure of
maximal spectral type) can easily have non-homogeneous  Lebesgue
spectrum. For Koopmans operators there is still no examples, since
all calculated examples of transformations with an absolutely
continuous component of finite multiplicity have a Lebesgue
component of the same multiplicity (see [1],[12],[13]).

The main reason why this definition is not so good is the following
theorem.

\begin{theorem}
There is no such transformations.
\end{theorem}

This implies that we really need to weaken the above conditions. In
our opinion, the weakest setting comes if we treat Lebesgue
components separately.

\begin{definition} We say that a transformation $T$ has
 \textbf{$k$-fold Lebesgue functions generated
spectrum } if for some collection of $k$-fold Lebesgue functions
$f_1,f_2,\ldots$
\[
L_0^2( \mu)=C(f_1,f_2,\ldots)
\]

\end{definition}

First note that, clearly, any $2$-fold Lebesgue functions generated
spectrum is just ordinary Lebesgue spectrum.

Now note that, if there exist  transformations with just one
$3$-fold Lebesgue function generated spectrum, then we can omit our
above discussion concerning relations between components for such
transformations. If so, it then can be considered as a solution of a
$k$-fold analog of the Banach problem in positive.

The main result is the following theorem.

\begin{theorem}
There is no  transformations with one $3$-fold Lebesgue function
generated spectrum.
\end{theorem}

Further note that, in spite of the case $k=2$, Examples 2.2 and 2.3
suggest that transformations with $3$-fold Lebesgue functions form a
bigger class. Indeed, recall one says  that a typical transformation
has some property (or the property is said to be typical) if the set
of elements satisfying the property contains a dense $G_{\delta}$
subset in the Polish group of all transformations equipped with weak
(coarse) topology. We offer the following simple theorem.
\begin{theorem}
For a typical transformation of the space $(X,\mathcal{F} , \mu)$,
there exist $k$-fold Lebesgue functions ($k=3,4,\ldots$). Moreover,
these functions can be chosen as finitely many-valued ones.
\end{theorem}

$Does $ $there$ $exist$ $a$ $transformation$ $with$ $no$ $multiple$
$Lebesgue$ $functions$ ?

We do not know.

Finally note that other certain questions of some fundamental origin
come almost automatically as well. To be precise, denote by $L(k)$
($L^*(k)$) the set of all transformations with $k$-fold Lebesgue
functions generated spectrum (with at least one $k$-fold Lebesgue
function)

$What$ $are $ $relations$ $between$  $L(k)$ or $L^*(k)$ $for$
$different$ $k$?

In particular,

$Is$ $it$ $true$ $that$ $L(2)\subset \cap_k L(k)$?

For bravery we can call these questions as Lebesgue analogs of
Rokhlin's problem. Working on this, we prove that $L(3)$ is not
included in $L(4)$ (see Example 2.4).

\section{  Key examples and proofs}

\begin{example} Bernoulli shifts.
\end{example}
Let $T$ be a Bernoulli shift. Take its standard realization on the
space
\[
(X,\mathcal{F} , \mu), \ X=\times_{i\in \mathbb{Z}}X_i, \
X_i=\mathbb{N}
\]
 for any $i$. Fix some orthogonal basis, say
$\{\varphi_i\}_i$, in the subspace of $L^2(X, \mu)$-functions
depending on the $0$-coordinate only. Clearly, we can assume that
every
 $\varphi_i$ is a bounded function and a constant function is one
 of the basis elements.

The reader will easily check that every non-constant function of the
form
\[
\widehat{T}^{i_1}\varphi_{j_1}\cdot\widehat{T}^{i_2}\varphi_{j_2}\cdots
\widehat{T}^{i_m}\varphi_{j_m} \ (-\infty<i_1<i_2<\ldots
<i_m<+\infty)
 \]
 is $k$-fold Lebesgue for any $k$, and,
consequently, every Bernoulli shift has $k$-fold Lebesgue functions
generated spectrum for any $k$.

It also should be noted that spectrum of any Bernoulli shift can not
be generated by finitely many  $k$-fold Lebesgue functions, since
its multiplicity function is unbounded.

\begin{example} Ergodic transformations with discrete spectrum.
\end{example}
Simple  $exotic$ examples of transformations with  Lebesgue and
mixing components appear if we do not require restrictions for
$k=2$. Indeed, take any ergodic transformation $T$ with discrete
spectrum. Then we get a certain  orthogonal decomposition
\[
L^2_0( \mu)=\bigoplus_{i\geq 1}C(f_i),
\]
where $\{f_1, f_2,\ldots \}$ is the set of all but the constant
eigenfunctions of $\widehat{T}$. Obviously, every $f_i$ is a
$3$-fold Lebesgue function and every $C(f_i)$ is a $3$-fold mixing
component. To get the same for any $k\geq3$, it is enough to
restrict ourselves to transformations with no roots of the unity
among eigenvalues of all their eigenfunctions. For example, every
ordinary ergodic shift on $\mathbb{R}/\mathbb{Z}$ is as it is
required.

\begin{example} of a   $3$-fold
mixing in $L^2_0( \mu)$ transformation that is not ($3$-fold)
mixing.
\end{example}
Let $T$ be any $3$-fold mixing transformation on the space
$(X,\mathcal{F} , \mu)$. Consider a
$\mathbb{Z}/2\mathbb{Z}$-extension, say $S$, of $T$ defined by
\[
S(x,y)=(Tx,y+1), \ (x,y)\in Y=X\times\mathbb{Z}/2\mathbb{Z};
\]
here $S$ preserves  the product measure, say $\nu$, of $\mu$ and the
Haar measure on $\mathbb{Z}/2\mathbb{Z}$. Clearly, $S$ is not
mixing, because there is a non-constant eigenfunction.

We claim that $S$ is $3$-fold mixing in $L^2_0( Y,\nu)$. Indeed,
first note that every function $\varphi\in L^2_0(Y, \nu)$ is the
orthogonal sum $\varphi_++\varphi_-$ according to the decomposition
$L^2_0(Y, \nu)=H_0^+\oplus H_0^-$ into $\widehat{S}$-invariant
components
\[
H_0^\pm=\{f\in L^2_0(Y, \nu): f(x,y+1)\equiv \pm f(x,y)\}.
\]
Clearly, $\widehat{S}\varphi_\pm=(\pm 1) \widehat{T^*}\varphi_\pm$;
 here the transformation $T^*$ is defined by $T^*(x,y)=(Tx,y)$.
 Take then any triple $f,g,h\in L^2_0(Y, \nu)$ of bounded functions. By
 the simple calculation,  for any  $n_1(m), n_2(m)$

($n_1(m)-n_2(m)\rightarrow\infty$ and $n_i(m)\rightarrow\infty$
$i=1,2$ as $m\rightarrow +\infty$ ) we get
\[
\int_Y\widehat{S}^{n_1(m)}f\cdot \widehat{S}^{n_2(m)}g \cdot
\overline{h} d\nu \rightarrow 0 \ \mbox{as } m\rightarrow +\infty,
  \]
since each of $8$ summands $\int_Y\widehat{S}^{n_1(m)}f_\pm\cdot
\widehat{S}^{n_2(m)}g_\pm \cdot \overline{h}_\pm d\nu$ is zero or
tends to zero as $m\rightarrow +\infty$.

\begin{example} of an ergodic transformation whose $3$-fold Lebesgue
functions generated spectrum  can not be $4$-fold Lebesgue functions
generated.
\end{example}
Suppose $\alpha$ is an irrational; then by $\Lambda$ we denote a
 countable subgroup of $\mathbb{T}=\{\lambda \in \mathbb{C} : |\lambda|=1\}$
 generated by $-1$ and $\exp (2\pi i\alpha)$. It is well known that
 there exists a unique (up to an isomorphism) ergodic
 transformation, say $T$, with discrete spectrum such that all its
 eigenvalues are $\Lambda(T)=\Lambda$. Arguing as in Example 2.2, we
 see that all but $f_1$ eigenfunctions $f_\lambda$ of $T$ are
$3$-fold Lebesgue, where $f_\lambda$ are eigenfunctions of
$\widehat{T}$ with eigenvalues $\lambda\in \Lambda$. Thus $T$ has
$3$-fold Lebesgue functions generated spectrum.

Let us prove that the spectrum of $T$ can not be $4$-fold Lebesgue
functions generated. To get this, it is sufficient to show that
every $4$-fold Lebesgue function is orthogonal to $f_{-1}$.

Assume the converse. Then $-1\in \mbox{supp } f_L$, where
$f_L=\sum_\lambda c_\lambda f_\lambda$, and $\mbox{supp }
f_L=\{\lambda\in \Lambda: c_\lambda\neq 0\}$. By the spectral
theorem, for any $\varphi, f \in L^2( \mu)$
\[
\varphi \in C(f) \Leftrightarrow \mbox{supp }\varphi\subseteq
\mbox{supp } f, \mbox{and }
\]
\[
C(\varphi) \perp C(f) \Leftrightarrow \mbox{supp }\varphi\cap
\mbox{supp } f=\emptyset.
\]
 Therefore $f_{-1}\in C(f_L)$.

 Besides, it is well known that every ergodic transformation $S$ with
  discrete spectrum is \textbf{rigid} (i.e.
  \[
\widehat{S}^{k_i}\rightarrow E \ \mbox{for some } k_i\rightarrow
\infty
\]
in the weak (strong) operator topology, where $E$ is the neutral
element).

It implies that for any pair
   of mutually different non-zero
integers ($n_2,n_3$)
\[
\forall \ n \int_X\widehat{T}^nf_L\cdot \widehat{T}^{n_2}f_L
   \cdot \widehat{T}^{n_3}f_L \cdot \overline{f}_L d\mu=0.
   \]
   Thus
   \[
   C(f_L)\bot C(\widehat{T}^{n_2}\overline{f}_L
   \cdot \widehat{T}^{n_3}\overline{f}_L \cdot f_L).
   \]
   It gives
\[
\forall \ n_2\neq n_3\neq 0\neq n_2  \ \int_Xf_{-1}\cdot
\widehat{T}^{n_2}f_L
   \cdot \widehat{T}^{n_3}f_L \cdot \overline{f}_L d\mu=0.
   \]
Applying rigidity again, we get
\[
   C(f_L)\bot C(\overline{f }_{-1}
   \cdot \widehat{T}^{n_3}\overline{f}_L \cdot f_L);
   \]
here $\overline{f }_{-1}
   \cdot \widehat{T}^{n_3}\overline{f}_L \cdot f_L
   \in L^2( \mu)$, since $|f_{-1}|$ is a constant function. Then

\[
\forall \ n_3\neq 0 \ \int_Xf_{-1}\cdot f_{-1}
   \cdot \widehat{T}^{n_3}f_L \cdot \overline{f}_L d\mu=0,
   \]
and, consequently,
\[
\int_Xf_{-1}\cdot f_{-1}
   \cdot f_{-1} \cdot \overline{f}_L d\mu=0, \
\int_X f_{-1} \cdot \overline{f}_L d\mu=0.
\]
This contradiction concludes the proof.

\begin{proof} $of $ $Proposition$ $1.4$. It is an easy application of
von Neumann ergodic theorem for unitary operators and left to the
reader.
\end{proof}

\begin{proof} $of $ $Theorem$ 1.8. It is an easy application of the
fact that a typical transformation $T$ is isomorphic to  a
$G$-extension for any finite abelian group $G$ (see [2]). Namely,
take a realization of $T$ as a $\mathbb{Z}/m\mathbb{Z}$-extension
for an
 $m$. Hence there is the standard decomposition of  $L^2( Y,\nu)$
into the orthogonal sum of $T$-invariant components
 \[
H_\chi=\{f \in L^2( Y,\nu): f(x,y)=\chi(y)\cdot \varphi(x),
\varphi\in  L^2( X,\mu) \},
\]
 where $\chi(y)$ are characters of
$\mathbb{Z}/m\mathbb{Z}$, $Y=X\times \mathbb{Z}/m\mathbb{Z}$. Every
non-zero element of $H_\chi$ ($\chi\neq 1$) is a $k$-fold Lebesgue
function for an appropriate $m$.
\end{proof}

\begin{proof}$of $ $Theorems$ 1.5 $and$ 1.7 . It is an easy application of the
following technical theorem
\begin{theorem}
Let $f$ be a $3$-fold Lebesgue function for a transformation $T$;
then there exists $n\neq 0$ such that $\widehat{T}^nf\cdot f\notin
C(f)\oplus C(1)$.
\end{theorem}
\end{proof}

The proof of Theorem 2.5 will be published elsewhere.

\section{Closing remarks and questions}

Looking at transformations treated above that are mostly of discrete
or Lebesgue maximal spectral types, one can think that dealing with
any $k$-fold Lebesgue function $f$ is just a matter of some power of
$f$ and convolutions for spectral types. However, applying group
extensions, one can easily  construct (see constructions in [4],
[7], [10]) transformations with
$\sigma_f\bot\sigma_{f^2}\bot\sigma_f \ast\sigma_f\bot\sigma_f$. It
also should be noted that for a typical transformation any measure
of the maximal spectral type is disjoint with its convolution. This
means that $\sigma_g\bot\sigma_f \ast\sigma_f$ for every pair
$f,g\in L_0^2( \mu)$.

 Transformations and, a bit more generally, group
actions (i.e.
 group representations by transformations) are main objects
  to study in modern
ergodic theory. Investigations describing roughly what is changed if
we go to actions of larger groups or back are one of the steadily
well-developing aspects of ergodic theory. The reader can easily
produce the notion being a ($k$-fold) Lebesgue function, as in the
case $k=2$, to actions of any countable groups. Traditional
($k$-fold) problems can be stated for general group actions as well.

Suppose a transformation $T$ belongs to $L(k)$; then by $m_k(T)$ we
 denote the minimal cardinality of collections of $k$-fold Lebesgue
 functions that generate the spectrum of $T$.

 $What$ $can$ $we$ $say$ $about$ $possible$ $values$ $of$ $m_k(T)$?

 Let us mention that, for example, $m_3(T)$ is infinite for every
 ergodic transformation $T$ with discrete spectrum. Related question
 is based on an easy inequality
 \[
m_2(T)\leq m_k(T),
\]
if they both exist.

$Can$ $we$ $say$ $the$ $same$ $for$ $other$ $pairs$ $n\leq k$?

It should be also noted that all above questions are most important
 in the main subclass of weakly mixing transformations.
It would be interesting to see any new effects that restriction may
imply.

\bibliographystyle{amsplain}

\end{document}